\documentclass[11pt]{article}
\usepackage[english]{babel}
\usepackage{latexsym, amssymb}
\usepackage{amsmath, color}
\usepackage{amsfonts}

\newtheorem{lemma}{Lemma}
\newtheorem{theorem}{Theorem}

\begin{document}
%\newenvironment{proof}{Proof}[section]
%\def\theequation{\arabic{section}.\arabic{equation}}

%-------------------------------------------------------------------
\def\a{\alpha}   \def\g{\gamma}  \def\b{\beta}  \def\th{\theta}
\def\l{\lambda}  \def\w{\omega}   \def\W{\Omega} \def\e{\varepsilon}
\def\d{\delta}   \def\f{\varphi}  \def\D{\Delta}    \def\r{\rho}
\def\s{\sigma}   \def\G{\Gamma}   \def\L{\Lambda}
\def\pd{\partial} \def\ex{\exists\,}  \def\all{\forall\,}
\def\dy{\dot y}   \def\dx{\dot x}   \def\t{\tau}
\def\bx{\bar x}   \def\bu{\bar u}   \def\bw{\bar w}  \def\bz{\bar z}
\def\p{\psi} \def\bh{\bar h}
\def\calO{{\cal O}}   \def\calP{{\cal P}}  \def\calR{{\cal R}}
\def\calN{{\cal N}}   \def\calL{{\cal L}}   \def\calQ{{\cal Q}}
\def\hx{\hat x}  \def\hw{\hat w}   \def\hu{\hat u}  \def\hv{\hat v}
\def\hz{\hat z}   % \def\ht{\hat t}

\def\empty{\O}     %\def\empty{\emptyset}
\def\lee{\;\le\;}   \def\gee{\;\ge\;}
\def\le{\leqslant}  \def\ge {\geqslant}
\def\Proof{{\bf Proof.}\quad}

\def\ov{\overline}  \def\ld{\ldots}
\def\tl{\tilde}   \def\wt{\widetilde}  \def\wh{\widehat}
\def\lan{\langle}  \def\ran{\rangle}

\def\hsp{\hspace*{20mm}} \def\vs{\vskip 8mm}
\def\bs{\bigskip}     \def\ms{\medskip}   \def\ssk{\smallskip}
\def\q{\quad}  \def\qq{\qquad}    \def\vad{\vadjust{\kern-5pt}}
\def\np{\newpage}   \def\noi{\noindent}
\def\pol{\frac12\,}    \def\const{\mbox{const}\,}
\def\R{\mathbb{R}}   \def\ctd{\hfill $\Box$}  \def\inter{\rm int\,}

\def\dis{\displaystyle}
\def\lra{\longrightarrow}  \def\iff{\Longleftrightarrow}
\def\beq{\begin{equation}\label}  \def\eeq{\end{equation}}
\def\bth{\begin{theorem}\label}   \def\eth{\end{theorem}}
\def\ble{\begin{lemma}\label}     \def\ele{\end{lemma}}
\def\begar{\begin{array}} \def\endar{\end{array}}

\newcommand{\vmin}{\mathop{\rm vraimin}}
\newcommand{\vmax}{\mathop{\rm vraimax}}

\newcommand{\blue}[1]{\textcolor{blue}{#1}}
\newcommand{\red}[1]{\textcolor{red}{#1}}
\newcommand{\gr}[1]{\textcolor{green}{#1}}

\def\bls{\baselineskip}   \def\nbls{\normalbaselineskip}
%\bls =1.1\nbls
%----------------------------------------------------------------

\begin{center}
{\Large\bf Reply to Comments on Our Paper \\ "On the Relation Between
Two Approaches to Necessary \\[3pt]
Optimality Conditions in Problems with State Constraints"
}
\\[8mm]
{\Large Andrei Dmitruk\footnotemark[1], \,\,  Ivan Samylovskiy\footnotemark[2]} \\[10pt]
%{\large А.В. Дмитрук\footnotemark[1], \quad Н.П. Осмоловский \footnotemark[2]}
%\\[20pt]
%{\large Moscow State University}\\[4pt]
\footnotemark[1]{Russian Academy of Sciences, CEMI, \\
Lomonosov Moscow State University, avdmi@cemi.rssi.ru}\\
\footnotemark[2]{Lomonosov Moscow State University, ivan.samylovskiy@cs.msu.ru}\\[20pt]
\end{center}

\begin{abstract}
In \cite{Kar}, a number of critical claims were made concerning our paper
\cite{DS}. \\ \hspace*{5mm} Here we reply to these claims.
\end{abstract}
\ms

%%---------------------------------------
\section*{Reply to Comments by Karamzin}

First of all, a highly biased style of these comments leaps to the eye.
%Setting aside the author's highly biased style, one can see that his claims
However, let us set it aside and concentrate on mathematics.
One can see that the author's claims can be essentially reduced to the following:
a) lack of novelty (that our results "can simply be derived from the already
known results in the literature"), and  b) incompleteness of the result
(stationarity conditions instead of Maximum Principle).
Let us make sense of this.
\ssk

As was clearly said in our paper \cite[Remark 9.1]{DS}, its novelty is a complete
realization of the idea of obtaining optimality conditions (in this case,
conditions of stationarity) in problems with a state constraint by differentiating
the state constraint on the contact interval and thus by passing to a problem
with a mixed constraint of the equality type.  This idea was first proposed
by Gamkrelidze \cite{Pont}, but in our paper it was realized, for a specific class
of problems, by another approach.
It needs to be emphasized that, rather than relying on  Gamkrelidze's
{\em result}, in order to obtain optimality conditions in the form of
Dubovitskii and Milyutin \cite{DM65} we realized {\em his idea}\, for obtaining
these conditions. (Perhaps, this should have been said in \cite{DS} more clearly.)
Our method is based on the replication of the state and control variables in
accordance with the number of qualitatively different subarcs of the reference
trajectory (interior or boundary w.r.t. the state constraint).
This replication trick is quite natural and has been known for a long time
(we give several references in \cite{DS}), but to our knowledge it has never
been used in the given situation with complete details and clarity
(and with a new feature consisting of a "two-stage varying" of the
reference trajectory). Despite the author's assertion, it is not the
"reduction to $v-$problem"\, in the sense of the $v-$problem proposed by
Dubovitskii and Milyutin \cite{DM65} and used in \cite[Sec. 2.5]{IT}.\,
The book \cite{IT} does not contain the replication method.
\ssk

Note that this replication trick turned out to be effective in some other
problems as well, such as problems with the control system of hybrid type,
with intermediate constraints, and with a variable structure,
where it allows to reduce these problems to a standard optimal control problem
and to apply the already known optimality conditions, but rather surprisingly,
it was missed even by highly prominent mathematicians, who obtained optimality
conditions in these problems by performing all the heavy procedure of variational
analysis (see comments and references in \cite{DK-SCL, DK-NCS}).
Therefore, in our opinion, this method is quite worthy of attention.
\ssk

Now, let us dwell on the paper by Neustadt \cite{Nst}. The author of the note
claims that its results allow one to make a passage from the conditions
of Gamkrelidze (GamC)  to those of Dubovitskii and Milyutin (DMC) by a change
of the costate variable, proposed in \cite{Nst}.  However, this is not true,
because the proposed passage contains a vicious circle. The matter is that
Neustadt's paper is not based either on the result of Gamkrelidze,
or on his idea of differentiating the state constraint on the contact
interval.  His approach is in fact close to that of Dubovitskii and Milyutin;\,
namely, just like they, he considers the state constraint as an inequality
constraint in the space of continuous function, and exactly from here, making
the linearization of all the constraints and the cost, and applying then the
Farkas-Minkowski theorem, he directly obtains a sign-definite measure possessing
all the required properties.  He clearly says (p. 134) that his conditions
{\em coincide}\,  with DMC and {\em imply}\,  GamC.
In fact, he makes a passage from DMC to GamC. (The same passage is later
performed by Arutyunov, Karamzin, and Pereira in \cite{AKP}, with no reference
to \cite{Nst}.)  But the reverse passage  GamC  $\to$ DMC  (by the reverse
transformation of the costate) is valid only if the measure in GamC is completely
sign-definite (and in this case one does not need to rely on paper \cite{Nst},
since the transformation, as well as the calculations in \cite{Kar}, is quite
simple), whereas GamC  do not provide the information about the sign of the
measure at its atoms.\,

Thus, in the claim concerning novelty Mr. Karamzin misleads the reader by
incorrect interpretations of methods and results in \cite{IT, Nst, AKP}.
%Those were our novelty concerns.
\ms

Now let's address the alleged incompleteness of our result. Note again that
we confined the study to the {\it stationarity conditions} (i.e. necessary first
order conditions) for the extended weak minimality.  Such conditions constitute
an important stage in the study of any optimization problem.  (Recall here
the Fermat condition  $f'(x_0) = 0$ in the problem  $f(x)\to \min,$
the Euler and Euler-Lagrange equations in the calculus of variations,
the Lagrange multipliers rule in the nonlinear programming, etc.)\,
Like any classical necessary first order conditions, they are not complete
in the sense that they can be strengthened  (e.g. by second order conditions).
But {\em as the stationarity conditions}\, they are complete, since they are
equivalent, loosely speaking, to the non-negativity of the cost derivative
in all admissible directions, and because of this their importance remains
undisturbed even if they are strengthened  by other conditions.
(The Weierstrass necessary condition %equivalent to the Maximum Principle,
does not invalidate the Euler equation!)\, % which is the stationarity condition!)\,
The issue of obtaining Maximum Principle {\em by the given approach}\,
was not considered in our paper, it requires additional study. \,
So, this claim of the author is also inadequate.  	
\ssk

Obviously, the presence of extended weak minimality implies the strong minimality
in an $\e$-tube around the reference control, and hence, the fulfillment
of Maximum Principle in this $\e$-tube, which in turn implies the absence
of atoms of the measure at the junction points with the state boundary.
In view of this, our sentence in \cite[p. 407]{DS}, given below in bold italics,
is, indeed, not correctly formulated:

{\it "Studies show ... that in case of strong (or at least Pontryagin type}
[16,17]) {\it minimality, the adjoint variable and measure do not have
jumps under condition} (2).
{\bf \em However, this result is not, in general, valid in the case of
extended weak minimality} {\it (the reason is that one cannot rely upon the
maximality of Pontryagin function w.r.t. u, having in disposal only the stationarity
of the extended Pontryagin function)."} \ssk

Here we should have said more precisely:\,  {\bf \em "However, this result is not,
in general, valid in the case of stationarity"}.  In fact, this is then said
in brackets.  Note also that in p. 410  we give a correct resume:\,

{\it "Thus, the stationarity conditions do not guarantee the absence of atoms,
while, according to} [5,15,19--23], {\it the maximum principle does."}\,
We hope that the attentive reader will be able to understand this point properly.\,
%The other claims in \cite{Kar} do not deserve attention.
\ssk

Having said that {\it "the technique of reducing to a mixed constrained problem is
obviously too restrictive as important information on the admissible trajectories
subject to state constraints is lost by this transition"}, Mr. Karamzin missed
the fact that this is just the first stage of our variation procedure. At the
second stage we use variations that go inside the state constraints.
\ssk

As concerns an omission in our assumptions, of course the data functions $f'_u$
and $g'_u$ should be assumed Lipschitz continuous, not just continuous.
\ssk

The author's claims concerning "wrong" citations and "confusing" title
present just his personal opinion.

%%-------------------------------------------------------

\end{document}